\documentclass[12pt]{article}
\usepackage{amssymb,amsmath,amsthm,psfig,amscd,xy,mathrsfs}

\def\C{\mathbb C}

\def\E{\mathbb E}

\def\I{\mathbb I}
\def\K{\mathbb K}
\def\L{\mathbb L}
\def\P{\mathbb P}
\def\Q{\mathbb Q}
\def\R{\mathbb R}

\def\Y{\mathbb Y}
\def\Z{\mathbb Z}

\def\Ge{\varepsilon}
\def\Gg{\gamma}
\def\GG{\Gamma}

\def\GD{\Delta}

\def\Go{\omega}

\def\GS{\Sigma}

\def\Gt{\tau}

\def\BE{\mathbf E}
\def\BF{\mathbf F}

\def\Bj{\mathbf j}
\def\Bs{\mathbf s}

\def\BT{\mathbf T}
\def\BV{\mathbf V}

\def\BGg{\boldsymbol\gamma}

\def\BGp{\boldsymbol\phi}

\def\BGt{\boldsymbol\tau}

\def\cA{\mathcal A}
\def\cB{\mathcal B}

\def\cI{\mathcal I}

\def\cR{\mathcal R}

\def\sR{\mathscr R}

\def\fM{\mathfrak M}
\def\fr{\mathfrak r}

\def\dd{\partial}
\def\Bdd{\boldsymbol\dd}
\def\smin{\setminus}

\def\emp{\emptyset}
\def\eksi{\mathbf{id}}

\def\im{\sqrt{-1}}

\newtheorem{thm}{Theorem}
\newtheorem{thmquote}{Theorem}

\newtheorem{lem}{Lemma}[section]

\newtheorem{cor}[thm]{Corollary}

\newtheorem{prop}[lem]{Proposition}

\theoremstyle{definition}

\newtheorem{exa}[lem]{Example}
\theoremstyle{remark}

\newtheorem{rem}[lem]{Remark}

\def\centerfig#1{\centerline{\psfig{figure=#1,silent=}}}

\newcommand{\ve}[1]{\mathbf{#1}}

\def\curve{(\GS;\ve{s})}

\newcommand{\cmod}[1]{\overline{M}_{#1}}
\newcommand{\cmodo}[1]{M_{#1}}

\newcommand{\konj}[1]{\overline{#1}}

\newcommand{\ove}[1]{\widetilde{#1}}

\newcommand{\cls}[1]{[\konj{D}_{#1}]}

\begin{document}

\title
{Chow groups of the moduli spaces of weighted pointed stable curves of genus zero}

\author{\"Ozg\"ur Ceyhan}



\maketitle

\begin{abstract}
The moduli space $\cmod{\cA}$ of weighted pointed stable curves of 
genus zero is stratified according to the degeneration types of such 
curves. We show that the homology groups of $\cmod{\cA}$ are generated 
by the strata of $\cmod{\cA}$ and give  all additive relations between them. 
We also observe that the Chow groups $A_i \cmod{\cA}$ and the homology 
groups $H_{2i} \cmod{\cA}$ are isomorphic. 
This generalizes Kontsevich-Manin's and Losev-Manin's theorems to 
arbitrary weight data $\cA$.
\end{abstract}



\section{Introduction}
\label{ch_intro}

The stability conditions have been used to compactify moduli spaces which have
the desired properties. Families of moduli spaces with respect 
to stability parameters have recently attracted attention in various 
contexts such as weighted maps \cite{ag,bm,mm1,mm2}, weighted pointed curves 
\cite{has,lm1,lm2,m1}, weighted configuration spaces, vector bundles 
and applications to representation theory \cite{kl1,kl2}, and  triangulated
categories \cite{b}.

Moduli space of pointed stable curves of genus $g$ has been extensively
studied in the literature. In \cite{has},  Hassett enriched the pointed 
curves by assigning a {\it weight} to each marked point and studied the
moduli stack of weighted pointed stable curves. In particular, he studied 
variations of the compactifications of moduli spaces as well as the 
corresponding chambers of stability conditions and wall crossing phenomena.

The homology groups of the moduli spaces $\cmod{\cA}$ of pointed weighted 
stable curves of genus zero corresponding to particular chambers of 
stability data play an essential role in the study of mirror symmetry, quantum
cohomology and Frobenius manifolds: The weight data $\cA= (1,1,\cdots,1)$
give the well known moduli space of pointed stable curves $\cmod{0,n}$ 
studied in \cite{ke,km1,km2,knu,m2}, and the weight data 
$\cA= (1,1,\Ge,\cdots,\Ge), 0 <\Ge \ll 1$ give the moduli space $L_{0,n}$ 
that has been studied by Losev and Manin in \cite{lm1,lm2,m1}.

In this note, we give a presentation of the homology groups of the moduli 
space $\cmod{\cA}$ for arbitrary weight data $\cA$. We show that the
homology groups are generated by the cycles of the strata of $\cmod{\cA}$,
as in the cases of Kontsevich-Manin and Losev-Manin. The additive
relations are obtained from the additive relations in the homology groups of
$\cmod{0,n}$ by using the reduction morphisms i.e., our description
generalizes Kontsevich-Manin's and Losev-Manin's theorems which give
additive structures on $H_* \cmod{0,n}$ and $H_* L_{0,n}$ to arbitrary 
weight data $\cA$. Since the homology groups are generated by the 
strata of $\cmod{\cA}$ and the relations are obtained from the relations
in $H_*(\cmod{0,n})$, the homology groups $H_{2i}(\cmod{\cA})$ are
isomorphic to the Chow groups $A_i(\cmod{\cA})$.  This result 
directly transfers the chamber decomposition and wall crossing phenomena 
to the homology level.

It is important to note that the technique used in this paper is significantly 
different than Keel's calculation of $A_* \cmod{0,n}$ in \cite{ke} and 
Musta\c{t}\u{a} \& Musta\c{t}\u{a}'s technique in
\cite{mm1,mm2}. Instead of using birational morphisms and blow-up
formulas for homology/Chow groups, we directly use the stratification of 
the moduli space $\cmod{\cA}$ and the spectral sequence of forgetful
morphisms. This technique has its own intrinsic power and may be
used beyond the cases where blow-up formulas are applicable; for instance,
it has been used in \cite{c}  to calculate the homology of the moduli space of
real curves where blow-up formulas can only calculate homology in 
$\Z/2\Z$  coefficients.


In this paper, all varieties are considered over the field $\C$ except when the 
contrary is stated. Therefore, we usually omit mentioning the base field. 

{\bf Plan of this paper:} In Section \ref{ch_w_curves}, we review some 
basic facts on weighted pointed curves of genus zero and their moduli 
space $\cmod{\cA}$. In the following section, we give the combinatorial 
stratification of $\cmod{\cA}$.
In Section \ref{ch_cohomology}, we consider  restrictions of 
forgetful morphisms to the strata of $\cmod{\cA}$ and study their 
fibers. We calculate the relative homology of the strata inductively
by using forgetful morphisms. Finally, in Section \ref{ch_hom_mod},
we give a combinatorial presentation of the homology and Chow groups 
of $\cmod{\cA}$ and prove the statement inductively by using the spectral 
sequence of the forgetful morphisms.

{\bf Acknowledgements:} I am thankful to Arend Bayer for valuable 
discussions and his careful reading. I am indebted to Sarah Carr and 
Andy Wand for their proofreading and commenting on my texts. I also wish to
thank to Anca Musta\c{t}\u{a} for bringing her papers to my attention.
Thanks are also due to
Max-Planck-Institut f\"ur Mathematik/Bonn and l'Universit\'e de 
Montr\'eal for their hospitality.

\section{Weighted pointed stable curves and their moduli}
\label{ch_w_curves}

This section reviews the basic facts on the moduli problem of weighted 
pointed stable curves of genus zero.

\subsection{The moduli problem}
\label{sec_n_curves}

A family of nodal curves of genus zero with $n$ {\it labeled points} over 
$B$ consists of 
\begin{itemize}
\item a flat proper morphism $\pi: C \to B$ whose geometric fibers $\GS$ 
are nodal connected curves of arithmetic genus zero; and%
\item a set of sections $\Bs=(s_1,\cdots,s_n)$ of $\pi$.
\end{itemize}

A {\it weight datum} $\cA$ is an element $(m_1,\cdots,m_n) \in \Q^n$ such that 
$0 < m_i \leq 1$ for $i =1,\cdots,n$, and $m_1+\cdots+m_n >2$.

A family of nodal curves of genus zero with $n$ labeled points $\pi: C \to B$ is 
{\it stable} (with respect to $\cA$) if 
\begin{itemize}
\item the sections $s_1,\cdots,s_n$ lie in the smooth locus of $\pi$, and 
for any subset $\{s_{i_1},\cdots,s_{i_r}\}$ with non-empty intersection we 
have $m_{i_1} + \cdots + m_{i_r} \leq 1$;
\item $K_{\pi}+ m_1 s_1 +\cdots + m_n s_n$ is $\pi$-relatively ample.
\end{itemize}

An  {\it $\cA$-pointed curve}  is a fiber $\curve$ of a family which is stable with 
respect to $\cA$.

\begin{thmquote}[Hassett, \cite{has}]
\label{thm_mod}
For any  $\cA$ with $n \geq 3$,  there exists a Deligne-Mumford
stack $\fM_{\cA}$, smooth and proper over $\Z$, representing the moduli 
problem of $\cA$-pointed curves. The corresponding moduli scheme 
$\cmod{\cA}$ is projective over $\Z$. 
\end{thmquote}

%
%

\subsection{Natural transformations}
\label{natural}

\paragraph{Reduction morphism:} 
Let $\cA=(m_1,\cdots,m_n)$ and $\cA'=(m'_1,\cdots,m'_n)$ be a pair 
of weight data such that $m'_i \leq m_i$ for all $i$.  In \cite{has}, Hassett 
showed that there exists a natural birational reduction morphism 
\begin{eqnarray}
\rho_{\cA,\cA'}: \cmod{\cA} \to \cmod{\cA'}.
\end{eqnarray}
The image $\rho_{\cA,\cA'} \curve$ of $\curve \in \cmod{\cA}$ is obtained by 
successively collapsing components of $\GS$ along which 
$K_\GS + m'_1 s_1 + \cdots +  m'_n s_n$ fails to be ample.

\paragraph{Forgetful morphism:} 

Let $\cA$ be a weight datum and $\cB= ((m_i)_{m_i \in \cA},m_{n+1})$.  
In \cite{has}, 
Hassett also showed that there exists a natural forgetful morphism 
\begin{eqnarray}
\pi_{\cB,\cA}: \cmod{\cB} \to \cmod{\cA}.
\end{eqnarray}
The image $\pi_{\cB,\cA} \curve$ of $\curve \in \cmod{\cB}$ is obtained by 
forgetting the labeled point $s_{n+1}$ and successively collapsing components of 
$\GS$ along which $K_\GS + m_1 s_1 + \cdots +  m_n s_n$ fails to be 
ample.

It is important to note that Hassett's result is in fact much more general;
he proved the same statements for arbitrary genus. 
%
%

\section{Stratification of the moduli space $\cmod{\cA}$}
\label{sec_strata}

In this section, we first introduce the combinatorial structures encoding
the degeneration types of $\cA$-pointed curves. Then, we give a 
stratification of the moduli space $\cmod{\cA}$ in terms of these
combinatorial data. 

\subsection{Combinatorial types of weighted pointed curves}

\subsubsection{Graphs}

A {\it graph} $\GG$ is a collection of finite sets of {\it vertices} 
$\BV_{\GG}$ and {\it flags (or half edges)} $\BF_{\GG}$ with a 
boundary map $\Bdd_{\GG}: \BF_{\GG} \to \BV_{\GG}$ and an 
involution $\Bj_{\GG}: \BF_{\GG} \to \BF_{\GG}$ ($\Bj_{\GG}^{2} =\eksi$). 
We call $\BE_{\GG} = \{(f_{1},f_{2}) \in \BF_{\GG}^{2} \mid f_{1}= j_{\GG} f_{2}\ \& \
f_{1} \not= f_{2} \}$ the set of {\it edges}, and 
$\BT_{\GG} = \{f \in \BF_{\GG} \mid f = j_{\GG} f\}$ the set of {\it tails}. 
For a vertex $v \in \BV_{\GG}$, let $\BF_{\GG}(v) =\Bdd^{-1}_{\GG}(v)$ 
and $|v| = |\BF_{\GG}(v)|$ be the {\it valency} of $v$.

A {\it weighted graph} is a graph $\GG$ endowed with a map 
$\cA_\GG: \BF_\GG \to \Q \cap ]0,1]$ such that $\cA_\GG(f) = 1$ for all 
flags that are part of an edge i.e., for which $\Bj_\GG (f) \ne f$.

We think of a graph $\GG$ in terms of  its {\it geometric 
realization} $||\GG||$: Consider the disjoint union of closed intervals 
$\bigsqcup_{f_i \in \BF_{\GG}} [0,1]  \times f_i$ and identify $(0, f_i)$ 
with $(0, f_j)$ if $\dd_\GG f_i = \dd_\GG f_j$, and identify $(t, f_i)$ 
with $(1-t, j_\GG f_i)$ for $t \in ]0,1[$ and $f_i \ne f_j$. The geometric 
realization of $\GG$ has a piecewise linear structure.

A {\it (weighted) tree}  is a (weighted) graph whose geometric 
realization is connected and simply-connected. 

\subsubsection{$\cA$-trees}

Let $\Gg$ be a weighted tree. An {\it $\fr$-structure} is a function on $\BV_\Gg$   
associating to each vertex $v \in \BV_\Gg$ an equivalence relation 
$\thickapprox$ on the set of tails $\BT_\Gg(v)$ such that 
$ \sum_{f_i \in [f]} \cA(f) \leq 1$ for each equivalence class  $[f]$.

An {\it $\cA$-tree} is a weighted tree $\Gg$ that carries an {\it $\fr$-structure} 
and  satisfies $\sum_{f \in \BF_\Gg(v)} \cA(f) >2$ for each $v \in \BV_{\Gg}$.

We denote $\cA$-trees by $(\Gg,\fr), (\Gt,\fr)$ or by bold Greek characters  
$\BGg, \BGt$. When it is necessary to indicate different $\fr$-structures on the same
weighted tree, we use indices in parentheses (e.g., $\BGg_{(i)}$).

The $\fr$-structure of an $\cA$-tree $\BGg$ determines the {\it weight structure}
 $$
 \cA_v = ((m_{[f]})_{[f] \in O_{\fr(v)}}, (1)_{f \in \BE_{\Gg}(v)}) \ \ 
 \text{for}\  \ \forall v \in \BV_\Gg
 $$
where $O_{\fr(v)}$ is the set of equivalence classes of $\fr(v)$, and 
$m_{[f]} = \sum_{f_i \in [f]} m_{f_i}$.

\subsubsection{Morphisms of graphs}

Let $\Gg$ and $\Gt$ be  trees with $n$ tails.  A {\it morphism} between these trees $\BGp: \Gg\to \Gt$ is a pair of maps $\BGp_{\BF}: \BF_{\Gt} \to \BF_{\Gg}$ and $\BGp_{\BV}: \BV_{\Gg} \to \BV_{\Gt}$  satisfying the following conditions:
\begin{itemize}
\item $\BGp_{\BF}$ is injective and $\BGp_{\BV}$ is surjective. %
\item The following diagram commutes
$$
\begin{CD}
\BF_{\Gg}  @> {\Bdd_{\Gg}} >>   \BV_{\Gg} \\
@A{\operatorname{\BGp}_\BF}AA      @VV{\operatorname{\BGp}_\BV}V \\
\BF_{\Gt}   @> {\Bdd_{\Gt}}>>   \BV_{\Gt}.
\end{CD}
$$
\item $\BGp_{\BF} \circ \Bj_{\Gt} = \Bj_{\Gg} \circ \BGp_{\BF}$.
\item $\BGp_{\BT}:={\BGp_{\BF} |}_{\BT}$ is a bijection.
\end{itemize}

Each morphism induces a piecewise linear map on geometric realizations.

An {\it isomorphism} $\BGp: \BGg \to \BGt$ is a morphism where $\BGp_{\BF}$ and $\BGp_{\BV}$ are bijections.


\subsubsection{Dual trees of weighted pointed curves}
Let $\curve$ be an $\cA$-pointed curve and $\eta: \hat{\GS}  \to \GS$ 
be its normalization. Let $(\hat{\GS}_v;\hat{\ve{s}}_v)$ be the following 
pointed curve: $\hat{\GS}_{v}$ is a component of  $\hat{\GS}$, and 
$\hat{\ve{s}}_v$ is the set of points consisting of the preimages of 
{\it special} (i.e.,  labeled and nodal) points on $\GS_{v} := \eta(\hat{\GS}_v)$. 
The points $\hat{\ve{s}}_v = (s_{f_{1}}, \cdots,  s_{f_{|v|}})$ on $\hat{\GS}_v$ 
are ordered by the elements  $f_*$ in the set $\{f_1,\cdots,f_{|v|}\}$.

The {\it dual tree}  of an  $\cA$-pointed curve $\curve$ is an 
$\cA$-tree $\BGg$ consisting of following data:
\begin{itemize}
\item $\BV_{\Gg}$ is the set of  components of $\hat{\GS}$.

\item $\BF_{\Gg}$ is the set consisting of the preimages of special points.

\item $\Bdd_{\Gg}: f \mapsto v$ if and only if $s_f \in \hat{\GS}_v$.

\item $\Bj_{\Gg}: f \mapsto f$ if and only if $s_f$ is a labeled point, and
$\Bj_{\Gg}: f_1 \mapsto f_2$ if and only if $s_{f_1} \in \hat{\GS}_{v_1}$ and
$s_{f_2} \in  \hat{\GS}_{v_2}$ are the preimages of  the nodal point
$\GS_{v_1} \cap \GS_{v_2}$.

\item $\cA_\Gg(f)=m_f$ if $f \in \BT_\Gg$, and $\cA_\Gg(f)=1$ if $f \in 
\BF_\Gg \smin \BT_\Gg$.

\item a pair of tails $f_1,f_2 \in \BT_\Gg(v)$ are equivalent if and only if 
$s_{f_1} = s_{f_2}$. 
\end{itemize}


\subsection{Combinatorics of degenerations} 
\label{sec_degen}

The degenerations of $\cA$-pointed curves are encoded by the morphisms of
$\cA$-trees as follows.

\subsubsection{Contractions of edges}
Let $\curve$ be an $\cA$-pointed curve and let  $\BGg = (\Gg,\fr)$ be its dual tree. 
Consider the deformation of a nodal point of  $\curve$. Such a deformation 
of $\curve$ gives a {\it contraction}  of an edge $\BGg \mapsto \BGg/e$: 
Let $e=(f_{e},f^{e}) \in \BE_{\Gg}$ be the edge corresponding to the nodal 
point and $\Bdd_{\Gg}(e) = \{v_{e},v^{e}\}$, and consider the equivalence 
relation $\thicksim$ on the set 
of vertices, defined by: $v \thicksim v$ for all $v \in \BV_{\Gg} \smin \{v_{e},v^{e}\}$, 
and  $v_{e} \thicksim v^{e}$. Then, there is a weighted tree $\Gg/e$ whose 
vertices are $\BV_{\Gg}/\thicksim$ and whose flags are $\BF_{\Gg} \smin 
\{f_e,f^e\}$.  The involution, boundary and weight maps of $\Gg/e$ are the 
restrictions of $\Bj_{\Gg}$, $\Bdd_{\Gg}$ and $\cA_\Gg$. The $\fr$-structure
is the same as before the contraction.

\subsubsection{Identifications of tails}
Let $C \to B$ be a family $\cA$-pointed curves whose dual tree is 
$\BGt=(\Gt,\fr_{old})$.  Let $f_{i_1},\cdots,f_{i_r}$ be the tails corresponding 
to a set of labeled points  supported by the same component $\GS_v$ 
with $m_{i_1} +\cdots + m_{i_r}\leq 1$. Consider the limits of this family where  
corresponding sections $s_{f_{i_1}},\cdots,s_{f_{i_r}}$ intersect. Such a 
degeneration gives a new $\fr$-structure $\fr_{new}$ on the set of tails of $\BGt$:
The new equivalence relation $\thickapprox_{new}$ is  given by using the
old one $\thickapprox_{old}$ as follows
\begin{itemize}
\item $f_i \thickapprox_{new} f_j$ iff $f_i \thickapprox_{old} f_j$
for each pair $f_i, f_j \in \BT_{\Gt}$, and
  
\item $f_{i_k} \thickapprox_{new} f_{i_l}$ for all $f_{i_k}, f_{i_l} \in \{f_{i_1},\cdots,f_{i_r}\}$. 
\end{itemize}


We use the notation $\BGg < \BGt$ to indicate that either $\BGt$ is obtained 
by contracting some edges of $\BGg$, or $\BGg$ is obtained by identifying a
set of tails of $\BGt$.

\subsection{Stratification of the moduli space $\cmod{\cA}$}
\label{stratacomplex}

The stratification of $\cmod{\cA}$ according to the degeneration types of 
its elements is a direct consequence of Hassett's theorem. The following 
statement is almost a tautology due to the definition of dual trees of 
$\cA$-pointed curves.


\begin{prop}
\label{thm_strata}
\begin{enumerate}
\item For any $\cA$-tree $\BGg$, there exists a quasi-projective subvariety 
$D_{\BGg} \subset \cmod{\cA}$ of codimension 
$|\BE_\Gg|+(n-\sum_{\BV_\Gg} |O_{\fr(v)}| )$ 
parameterizing  $\cA$-pointed curves whose dual tree is $\BGg$. The subvariety 
$D_{\BGg}$ is isomorphic to $\prod_{v \in \BV_{\Gg}}  \cmodo{\cA_{v}}$ where 
$\cA_{{v}}$ is the weight structure at vertex $v$.  

\item $\cmod{\cA}$ is stratified by pairwise disjoint subvarieties $D_{\BGg}$. 
The closure $\overline{D}_{\BGg}$ of any stratum $D_{\BGg}$  is stratified by 
$\{D_{\BGg'} \mid \BGg' \leq \BGg \}$.
\end{enumerate}
\end{prop}

\begin{exa} 
Consider the case $n=4$. For any $\cA$ with $n=4$, the 
principal stratum $\cmodo{\cA}$ is isomorphic to $\P^1 \setminus \{0,1,\infty\}$. 
Moreover, its compactification $\cmod{\cA}$ is isomorphic to $\P^1$ obtained
by adding three boundary divisors. However, 
the stratification is encoded by different sets of $\cA$-trees. 

Figure \ref{example} depicts the $\cA$-trees used in the combinatorial stratification
of the  three separate cases: (a) $\cA = (1,1,1,1)$ i.e., $\cmod{\cA} = \cmod{0,4}$; 
(b) $\cA = (1,1,\Ge,\Ge)$ where $0< \Ge \ll 1$ i.e., $\cmod{\cA} = L_{0,4}$; 
(c) $\cA = (1,\Ge,\Ge,\Ge)$ where $ 1/3 < \Ge \leq 1/2$.
\begin{figure}[htb]
\centerfig{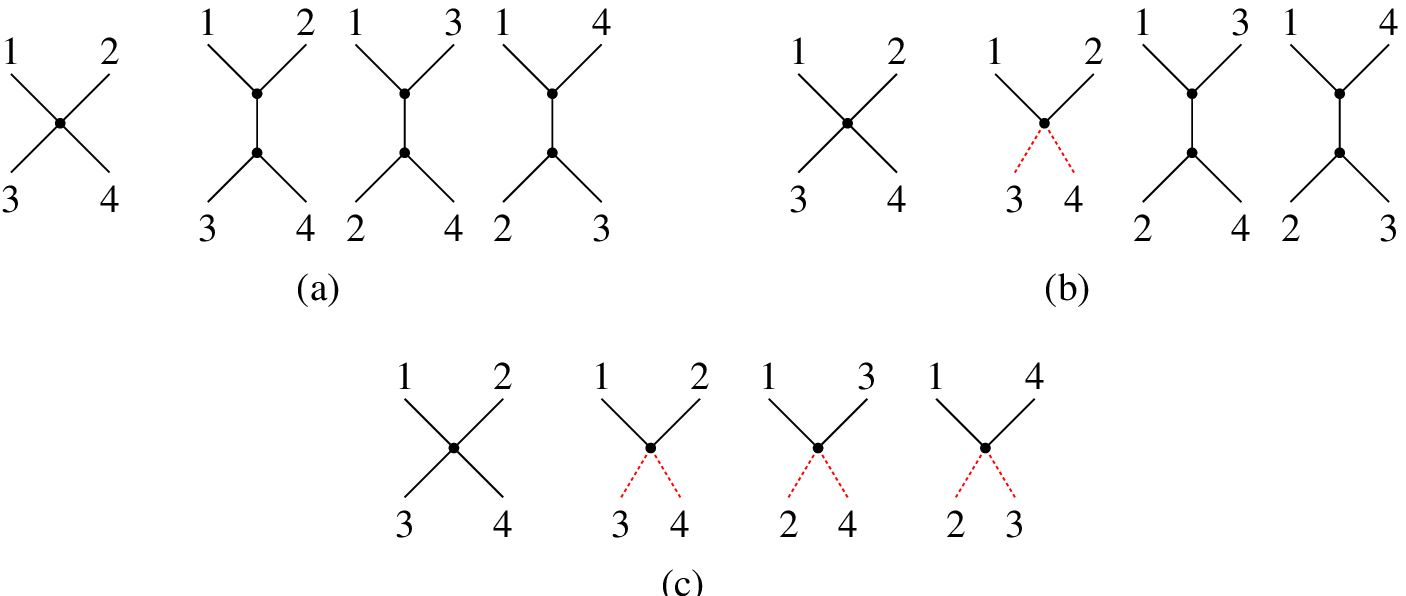,height=2.1in} %
\caption{(a) $\cA = (1,1,1,1)$;  
(b) $\cA = (1,1,\Ge,\Ge)$; 
(c) $\cA = (1,\Ge,\Ge,\Ge)$.}%
\label{example}
\end{figure}
In this figure, the dotted red tails depict the tails lying in the same equivalence class. 
\end{exa}

\section{Homology of the strata of $\cmod{\cA}$}
\label{ch_cohomology} 

In this chapter, we calculate the homology of the strata of $\cmod{\cA}$ relative to the union their substrata of codimension one and higher.




\subsection{Fibers of the forgetful morphism}
\label{sec_forget}

Let $\cA$ be a weight datum and $\cB= ((m_i)_{m_i \in \cA},m_{n+1})$. 
Let $\pi: D_{\BGg^*} \to D_{\BGg}$   be the restriction of the   
morphism $\pi_{\cB,\cA}: \cmod{\cB} \to \cmod{\cA}$ which forgets the labeled 
point $s_{n+1}$. Let $\BGp^{s_{n+1}}:\Gg^* \to\Gg$ be the corresponding 
forgetful morphism of trees, and let  $v_{s} = \Bdd_{\Gg^*}(s_{n+1})$.
In order to avoid the trivial cases here, we assume that $s_{n+1} \ne s_f$
for all $f \in \BT_{\Gg^*}(v_s)$.

We will denote the fibers $\pi^{-1} \curve$ of the forgetful morphism $\pi$ 
simply by $F_{\BGg^*}$.

\begin{lem} 
\label{lem_fil_I}
Let $\curve \in D_{\BGg}$.  Then, the fiber $F_{\BGg^*}$  is 
\begin{enumerate}
\item a projective line $\P^1$ minus the special points $s_{f}$ where $f \in \BF_{\Gg^*}(v_{s}) \smin \{s_{n+1}\}$  
if $\sum_{f \in \BF_{\Gg^*}(v_{s}) \smin \{s_{n+1}\}} \cB(f) > 2$;

\item the product $\prod_{v \in \BV_{\Gg^*} \smin (\BGp^{s_{n+1}}_\BV)^{-1} (\BV_{\Gg})} \cmodo{\cB_{v}}$ if 
$\sum_{f \in \BF_{\Gg^*}(v_{s}) \smin \{s_{n+1}\}} \cB(f) \leq 2$.
\end{enumerate}
\end{lem}

\begin{proof} 
Pick a curve $\curve \in D_{\BGg}$.  Let $(\GS^*,\ve{s}^*)$ be points in the fiber $F_{\BGg^*}$.

(1)
If $(\GS^*,\ve{s}^*) \in D_{\BGg^*}$ does not require the contraction of its component $\GS^*_{v_{s}}$ after 
forgetting the labeled point $s_{n+1}$, then the special points on $\GS^*_{v_{s}}$ must satisfy the inequality 
$\sum_{f \in \BF_{\Gg^*}(v_{s}) \smin \{s_{n+1}\}} \cB(f) > 2$. For all such $(\GS^*,\ve{s}^*)$, the curve $\GS^*$ is  
isomorphic to $\GS$, so $\GS^*_{v_{s}}$ (and $\GS_{v_{s}}$) is isomorphic to $\P^1$. All special 
points are  fixed (by choosing $\curve$) except $s_{n+1}$. The labeled point $s_{n+1}$ can be any point in 
$\GS_{v_{s}} \smin  \{s_f \mid f \in \BF_{\Gg^*}(v_{s}) \smin \{s_{n+1}\}\}$ since we assume that
$s_f \ne s_{n+1}$  for $f \in \BF_{\Gg^*}(v_{s}) \smin \{s_{n+1}\}$. Hence, the  elements $(\GS^*,\ve{s}^*)$ of 
$F_{\BGg^*}$ are  determined only by the position of the labeled  point  $s_{n+1}$ in $\GS_{v_{s}}$
i.e., the fiber is a punctured  $\P^1$ as stated above.

(2)
If   $(\GS^*,\ve{s}^*) \in D_{\BGg^*}$ requires the contraction of $\GS^*_{v}$ after forgetting $s_{n+1}$
for $v \in \BV_{\Gg^*} \smin (\BGp^{s_{n+1}}_\BV)^{-1} (\BV_{\Gg})$, then the special points supported by $\GS^*_{v}$ can realize all possible configurations. Therefore, the fiber $F_{\BGg^*}$ is $\prod_v \cmodo{\cB_{v}}$ where the product runs over all
contracted components.
\end{proof}

\subsection{Homology of the two dimensional fibers of the forgetful morphisms}
\label{sec_fib_hom}
Let $\pi: D_{\BGg^*} \to D_{\BGg}$ be the map forgetting the labeled point  $s_{n+1}$ 
which is discussed above.  Assume that $\sum_{f \in \BF_{\Gg^*}(v_{s}) \smin \{s_{n+1}\}} \cA(f) > 2$ 
i.e.,  the fibers are  projective lines  with punctures 
(see case (1) of  Lemma \ref{lem_fil_I}).

Let us identify the  fiber $F_{\BGg^*}$ with $\P^1 \smin \Bs$, and let 
$z=[z:0]$  be an affine coordinate on it. Let $z_f$ denote the positions of special 
points $s_f$ in these coordinates.

Then, the cohomology of a fiber is generated by the logarithmic differentials
\begin{eqnarray*}
H^0(F_{\BGg^*}) = \Z,  \ \ \ \
H^1(F_{\BGg^*}) = 
\bigoplus_{f}\ \Z \ \Go_{f ,s_{n+1} }
\end{eqnarray*}
where
\begin{eqnarray}
\label{eqn_logform}
\Go_{f ,s_{n+1} } = \frac{1}{2 \pi \im}\ d \log (z-z_f) 
\end{eqnarray}
for $f \in \BF_{\Gg^*} (v_{s}) \smin \{s_{n+1}\}$.

The homology with closed support $H^c_{1}(F_{\BGg^*})$ is isomorphic to the cohomology 
group $H^1(F_{\BGg^*})$ via Poincar\'e duality. Moreover, the definition of homology with 
closed support implies that $H^c_{1}(F_{\BGg^*})$  is isomorphic to the relative homology group 
$H_{1}(\overline{F}_{\BGg^*},\Bs)$ of the closure of the fiber. The group  
$H_{1}(\overline{F}_{\BGg^*},\Bs)$ is clearly generated by the homotopy classes of arcs 
connecting the pairs of punctures $s_{f_1}, s_{f_2}$. These arcs are the duals of the 
cohomology classes  $\Go_{f_1,s_{n+1}} - \Go_{f_2,s_{n+1}}$. We denote them by 
$\sR_{s_{n+1}, f_1 f_2}$. The homology groups $H^c_{2}(F_{\BGg^*})$  and 
$H_{2}(\overline{F}_{\BGg^*},\Bs)$  are isomorphic to  $H^0(F_{\BGg^*}) = \Z$. Hence
\begin{eqnarray*}
H_2(\overline{F}_{\BGg^*},\Bs) = \Z,  \ \ \ \ 
H_1(\overline{F}_{\BGg^*},\Bs) =
\left( \bigoplus_{f}\ \Z \ \sR_{s_{n+1},f_1 f_2}\right) / \cI_{\BGg^*}
\end{eqnarray*}
where  $f_i \in \BF_{\Gg^*}(v_{s}) \smin \{s_{n+1}\}$, and the subgroup of relations $\cI_{\BGg^*}$ is generated by 
\begin{eqnarray}
\label{eqn_f_rel}
\sR_{s_{n+1},f_1 f_2} + \sR_{s_{n+1},f_2 f_3} + \sR_{s_{n+1},f_3 f_1}.
\end{eqnarray}

\subsection{Homology of the strata}

In this section, we give the generators of the homology of the closed strata 
$\konj{D}_{\BGg}$  relative to the union of their lower dimensional strata 
$Q_{\BGg} := \bigcup_{\BGt < \BGg} D_{\BGt}$.

\begin{lem}
\label{lem_closed_hom}
Let $\pi: D_{\BGg^*} \to D_{\BGg}$ be the forgetful morphism discussed 
Section \ref{sec_forget}. Then, 
\begin{eqnarray*}
H^c_{d} (D_{\BGg^*};\Z) = \bigoplus_{p+q=d}
H^c_{p} (D_{\BGg};\Z) \otimes H^c_{q} (F_{\BGg^*};\Z).
\end{eqnarray*}
\end{lem}

\begin{proof} 
The strata $D_{\BGg^*}$ and $D_{\BGg}$ are given by the products
\begin{eqnarray*}
\prod_{v\in \BV_{\Gg^*}} \cmodo{\cB_{v}},\ \ \ \
\prod_{v\in \BV_{\Gg}} \cmodo{\cA_{v}}
\end{eqnarray*}
due to Proposition \ref{thm_strata}. The forgetful map $\pi$ preserves the components 
$(\GS^*_v,\ve{p}^*_v)$ of $(\GS^*,\ve{p}^*) \in D_{\BGg^*}$ except the contracted ones. 
Hence,  it gives the identity map on the factors $\cmodo{\cB_{v}}  \to \cmodo{\cA_{v}}$
for  $v \in (\BGp^{s_{n+1}}_\BV)^{-1} (\BV_{\Gg})$. On the other hand, it gives a fibration
\begin{eqnarray}
\label{eqn_fib_int1}
\pi_{res}: \prod_{v  \in \BV_{\Gg^*} \smin (\BGp^{s_{n+1}}_\BV)^{-1} (\BV_{\Gg})} \cmodo{\cB_{v_s}}   \to  point
&\mathrm{when}&  \sum_{f \in \BF_{\Gg^*}(v_{s}) \smin \{s_{n+1}\}} \cB(f) \leq 2, \\
\label{eqn_fib_int2}
\pi_{res}: \cmodo{\cB_{v_{s}}}   \to \cmodo{\cA_{v_{s}}}
&\mathrm{when} & \sum_{f \in \BF_{\Gg^*}(v_{s}) \smin \{s_{n+1}\}} \cB(f) > 2, 
\end{eqnarray}
with the same fibers $F_{\BGg^*}$ of $\pi: D_{\BGg^*} \to D_{\BGg}$.

In the case of $ \sum_{f \in \BF_{\Gg^*}(v_{s}) \smin \{s_{n+1}\}} \cB(f) \leq 2$, the stratum 
$D_{\BGg^*}$ is clearly $D_{\BGg} \times F_{\BGg^*}$. Hence, the claim follows from the 
K\"unneth formula.

In the case of $\sum_{f \in \BF_{\Gg^*}(v_{s}) \smin \{s_{n+1}\}} \cB(f) > 2$, 
the spaces $\cmodo{\cB_{v_{s}}}$ and $\cmodo{\cA_{v_{s}}}$ are diffeomorphic 
to the products of $\P^1 \smin \{s_1,s_2,s_3\}$ minus all diagonals. The map 
$\pi_{res}$ forgets the coordinate subspace $\P^1$ corresponding to the labeled 
point $s_{n+1}$ i.e., it is
\begin{eqnarray*}
(\P^1 \smin \{s_1,s_2,s_3\})^{|\BF_{\Gg^*}(v_{s})|-3} \smin \GD^*  
\to
(\P^1 \smin \{s_1,s_2,s_3\})^{|\BF_{\Gg}(v_{s})|-3} \smin \GD.
\end{eqnarray*}

The logarithmic differentials $d \log(z-z_f)$ give global cohomology classes 
on $\cmodo{\cB_{v_{s}}}$. On the other hand, we have seen that the restriction 
of these logarithmic forms to each fiber  generate the cohomology of the fiber 
(see Section \ref{sec_fib_hom}). By using  the Leray-Hirsch theorem, we obtain
\begin{eqnarray*}
H^{d} (D_{\BGg^*}) = \bigoplus_{p+q=d} H^{p} (D_{\BGg}) \otimes H^{q} (F_{\BGg^*}).
\end{eqnarray*}

The duality between cohomology  and homology with closed supports gives us 
isomorphisms which are needed to complete the proof.
\end{proof}

Since the strata of $\cmod{\cA}$ are the products given   in Proposition 
\ref{thm_strata}, their homology groups are obtained from the homology of their factors 
$\cmod{\cA_v}$ by using K\"unneth formula.  Here, we give only the relative 
homology for the strata of  one-vertex trees:

Now, let $\BGg$ be a one-vertex $\cA$-tree, and let $Q_{\BGg}$ be the  union of
 the codimension one and higher strata of the closed stratum $\konj{D}_{\BGg}$.

\begin{prop}
\label{lem_hom_strata1}
The relative homology group  
$H_{\dim (D_{\BGg}) -d} (\konj{D}_{\BGg},Q_{\BGg};\Z)$ is
generated by
\begin{eqnarray*}
\sR_{s_{i_1}, s_{j_1} s_{k_1}} \otimes \cdots \otimes \sR_{s_{i_d}, s_{j_d} s_{k_d}}
\end{eqnarray*}
where $j_*,k_* < i_*$ and $i_1 < \cdots < i_d  <n$. In particular,
\begin{eqnarray*}
H_{\dim (D_{\BGg})} (\konj{D}_{\BGg},Q_{\BGg};\Z) = \Z \ \cls{\BGg}
\end{eqnarray*}
where $\cls{\BGg}$ is the fundamental class of the stratum $D_{\BGg}$.
\end{prop}

\begin{proof}
The homology with closed support is defined by
\begin{eqnarray*}
H^c_*(D_{\BGg}) = \lim_{\to} H_*(D_{\BGg}, D_{\BGg}\smin K)
\end{eqnarray*}
where $K$ ranges over all closed subsets of  $D_{\BGg}$. The group 
$H_*(D_{\BGg},D_{\BGg} \smin K)$ is isomorphic 
to $H_*(\konj{D}_{\BGg}, \konj{D}_{\BGg}\smin \konj{K})$ where 
$\konj{K}$ ranges over all closed subsets of  $\konj{D}_{\BGg}$ which 
do not intersect  $Q_{\BGg}$. In the limit, 
$\konj{D}_{\BGg} \smin \overline{K}$ gives the union of substrata
$Q_{\BGg}$. Hence, the homology with closed support is isomorphic to 
the relative homology of $\konj{D}_{\BGg}$.

On the other hand, Lemma \ref{lem_closed_hom} implies that 
\begin{eqnarray*}
H_{d} (\konj{D}_{\BGg^*}, Q_{\BGg^*};\Z) = \bigoplus_{p+q=d}
H_{p} (\konj{D}_{\BGg}, Q_{\BGg};\Z) \otimes H_{q} (\overline{F}_{\BGg^*},\Bs;\Z). 
\end{eqnarray*}
where $F_{\BGg^*}$ is the fiber of the map $\pi: D_{\BGg^*} \to D_{\BGg}$ 
which forgets $s_{n+1}$.

We obtain the result by applying  the forgetful morphism successively 
and using the generators of  the relative homologies the 
fibers given in Section \ref{sec_fib_hom}.  In order to simplify the notation, 
we omit the factors coming from the generators of the second homology 
of the fibers.

It is clear that the top dimensional relative homology is generated by the relative 
fundamental class $\cls{\BGg}$.
\end{proof}

\section{Homology groups of $\cmod{\cA}$}
\label{ch_hom_mod}

In this section, we  give the homology groups of the moduli space 
$\cmod{\cA}$ in terms of generators and relations.

\subsection{Relations between strata of $\cmod{\cA}$}
\label{sec_ideal}
Here, we introduce a set of interesting relations between the strata of 
$\cmod{\cA}$ which will play a crucial role in our description of the homology 
of $\cmod{\cA}$.

Consider an $\cA$-tree $\BGg$ such that $\dim D_{\BGg}=d+1$, and a vertex  
$v \in \BV_\Gg$ with $|O_{\fr(v)}| \geq 4$. Let  $f_1,f_2,f_3,f_4 \in \BF_\Gg(v)$ and 
let their equivalence classes be pairwise disjoint. 
Put $\BF = \BF_\Gg(v) \smin \{f \mid f \thickapprox f_i, i=1,\dots 4\}$ and 
let $\BF_1,\BF_2$ be a partition of $\BF$ such that $\BF_i \bigcap [f]$ is 
either empty or equal to $[f]$ for all $f \in \BF$.

We define two $\cA$-trees $\BGg_1,\BGg_2$ such that $dim\ D_{\BGg_i} = d$:

\paragraph{The $\cA$-tree $\BGg_1$:} In order to define $\BGg_1$, we first
introduce a weighted tree $\widehat{\Gg_1}$. It is obtained by inserting a 
new edge $e=(f_{e},f^{e})$ into $\BGg$ at $v$ and its flags are given as follows:
Let $\Bdd_{\widehat{\Gg_1}}(e) =\{v_e,v^e\}$. The distribution of flags is 
given by $\BF_{\widehat{\Gg_1}}(v_e) = \BF_1 \cup [f_1] \cup [f_2] \cup \{f_e\}$ and 
$\BF_{\widehat{\Gg_1}}(v^e) = \BF_2 \cup [f_3] \cup [f_4] \cup \{f^e\}$.

Consider the following $\fr$-structure  on $\widehat{\Gg_1}$: Let the
equivalence relations on $\BT_{\widehat{\Gg_1}}(v_e)$, $\BT_{\widehat{\Gg_1}}(v^e)$ 
be the restrictions of the 
equivalence relation on $\BT_{\Gg}(v)$, and the equivalence relations 
for $\BT_{\widehat{\Gg_1}}(v) =\BT_{\Gg}(v)$ for  $v \in \BV_{\widehat{\Gg_1}} \smin \{v_e,v^e\}$ 
remain the same as before. Denote this $\fr$-structure by $\fr_1$. Then, we 
define the $\cA$-tree $\BGg_1$ in two separate cases:

\begin{itemize}
\item {\bf Stable case:} If $\widehat{\Gg_1}$ and the $\fr$-structure
$\fr_1$ provide an $\cA$-tree (i.e., $\sum_{f \in \BF_{\widehat{\Gg_1}}(v)} \cA(f) >2$ for  
$v \in \{v_e,v^e\}$), then we put $\BGg_1 := (\widehat{\Gg_1},\fr_1)$.

\item {\bf Unstable case:} If $\widehat{\Gg_1}$ and the $\fr$-structure
$\fr_1$ do not provide an $\cA$-tree (i.e., $\sum_{f \in \BF_{\widehat{\Gg_1}}(v_e)} \cA(f) \leq 2$ or  
$\sum_{f \in \BF_{\widehat{\Gg_1}}(v^e)} \cA(f) \leq 2$), then $\BGg_1$ is obtained from $\BGg$ 
by identifying the tails in $[f_1],[f_2]$ or $[f_3],[f_4]$ whichever are
adjacent to the unstable vertex (i.e., the one satisfying the weight inequality
$\sum_{f \in \BF_{\widehat{\Gg_1}}(v)} \cA(f) \leq 2$) in $\widehat{\Gg_1}$.
\end{itemize}

\paragraph{The $\cA$-tree $\BGg_2$:} The $\cA$-tree $\BGg_2$ is obtained
the same as $\BGg_1$ after swapping $f_2$ and $f_3$: We first
introduce a weighted tree $\widehat{\Gg_2}$. It is obtained by inserting a 
new edge $e=(f_{e},f^{e})$ into $\BGg$ at $v$ and its flags are given as follows:
Let $\Bdd_{\widehat{\Gg_2}}(e) =\{v_e,v^e\}$. The distribution of flags is 
given by $\BF_{\widehat{\Gg_2}}(v_e) = \BF_1 \cup [f_1] \cup [f_3] \cup \{f_e\}$ and 
$\BF_{\widehat{\Gg_2}}(v^e) = \BF_2 \cup [f_2] \cup [f_4] \cup \{f^e\}$.

Consider the following $\fr$-structure  on $\widehat{\Gg_2}$: Let the
equivalence relations on $\BT_{\widehat{\Gg_2}}(v_e)$, $\BT_{\widehat{\Gg_2}}(v^e)$ 
be the restrictions of the 
equivalence relation on $\BT_{\Gg}(v)$, and the equivalence relations 
for $\BT_{\widehat{\Gg_2}}(v) =\BT_{\Gg}(v)$ for  $v \in \BV_{\widehat{\Gg_2}} \smin \{v_e,v^e\}$ 
remain the same as before. Denote this $\fr$-structure by $\fr_2$. Then, we 
define the $\cA$-tree $\BGg_2$ in two separate cases:

\begin{itemize}
\item {\bf Stable case:} If $\widehat{\Gg_2}$ and the $\fr$-structure
$\fr_2$ provide an $\cA$-tree (i.e., $\sum_{f \in \BF_{\widehat{\Gg_2}}(v)} \cA(f) >2$ for  
$v \in \{v_e,v^e\}$), then we put $\BGg_2 := (\widehat{\Gg_2},\fr_2)$.

\item {\bf Unstable case:} If $\widehat{\Gg_2}$ and the $\fr$-structure
$\fr_2$ do not provide an $\cA$-tree (i.e., $\sum_{f \in \BF_{\widehat{\Gg_2}}(v_e)} \cA(f) \leq 2$ 
or  $\sum_{f \in \BF_{\widehat{\Gg_2}}(v^e)} \cA(f) \leq 2$), then $\BGg_2$ is obtained from $\BGg$ 
by identifying the tails in $[f_1],[f_3]$ or $[f_2],[f_4]$ whichever are
adjacent to the unstable vertex  in $\widehat{\Gg_2}$.
\end{itemize}

\paragraph{Principal relations between strata:}
We define $d$-dimensional class 
\begin{eqnarray} 
\label{eqn_rel}
\cR(\BGg;v,f_1,f_2,f_3,f_4) :=    
\sum_{\BGg_1}  \cls{\BGg_1} -  \sum_{\BGg_2}  \cls{\BGg_2}
\end{eqnarray}
where sum is taken over the isomorphism classes of all possible $\BGg_1$'s and $\BGg_2$'s 
for a fixed set of flags $\{f_1,f_2,f_3,f_4\}$.

\begin{lem}
\label{lem_relation}
$\cR(\BGg;v,f_1,f_2,f_3,f_4)$ is rationally equivalent to zero.
\end{lem}

\begin{rem}
If we consider only the fully stable cases (i.e., $\cA(f_1) + \cA(f_2) >1$, 
$\cA(f_3) + \cA(f_4) >1$, $\cA(f_1) + \cA(f_3) >1$ and $\cA(f_2) + \cA(f_3) >1$),
such as the cases where $\cA_\I= (1,1,\cdots,1)$  or where all of the $f_i$ are 
parts of edges, the 
homology relations given in Lemma \ref{lem_relation} reduce to the additive
relations in $H_*(\cmod{0,n})$ given by Kontsevich and Manin in \cite{km2}. They appear as 
a consequence of Keel's work \cite{ke}. In \cite{km2}, Kontsevich and Manin 
proved  that they in fact generate all additive relations in $H_*(\cmod{0,n})$.
\end{rem}

\begin{proof} (of Lemma \ref{lem_relation})
All additive relations in the Chow group $A_d(\cmod{0,n})$ are generated by 
\begin{eqnarray} 
\label{eqn_rel_n}
R(\BGg;v,f_1,f_2,f_3,f_4) :=    
\sum_{f_1 f_2 \ove{\BGg}_1 f_3 f_4}  \cls{\ove{\BGg}_1} -  
\sum_{f_1 f_3 \ove{\BGg}_2 f_2 f_4}  \cls{\ove{\BGg}_2} \equiv 0
\end{eqnarray}
where $\ove{\BGg}_i, i=1,2$ are  $\cA_\I$-trees and where $\cA_\I= (1,1,\cdots,1)$ 
(see, for example \cite{m2}). The notation $f_1 f_2 \ove{\BGg} f_3 f_4$ indicates  
the flags $f_1 f_2$ and $f_3 f_4$ are adjacent to two vertices $\{v_e,v^e\} = \dd_{\BGg}(e)$ of 
$\ove{\BGg}$.

Consider the reduction morphism  $\rho_{\cA_\I,\cA_{v}}: \cmod{0,|\BF_\Gg(v)|} \to \cmod{\cA_{v}}$ 
for the vertex $v$ of the $\cA$-tree $\BGg$. The push-forward  
$(\rho_{\cA_\I,\cA_{\Gg_v}})_* R(\BGg;v,f_1,f_2,f_3,f_4)$ gives us the relations 
$\cR(\BGg;v,f_1,f_2,f_3,f_4) \equiv 0$.
\end{proof}

\subsection{Homology groups of $\cmod{\cA}$}

\begin{thm}
The homology groups of $\cmod{\cA}$ are
\begin{eqnarray*}
H_{2d} (\cmod{\cA})   &=& 
\left( \bigoplus_{\BGg \mid \ \dim_\R D_{\BGg} =2d}  
H_{2d} (\konj{D}_{\BGg},Q_{\BGg}) \right) / \cI_d, \\
H_{2d-1} (\cmod{\cA})   &=& 0
\end{eqnarray*}
where the subgroup of relations $\cI_d$ is generated by $\cR(\Gg;v,f_1,f_2,f_3,f_4)$, see 
(\ref{eqn_rel}). 
\end{thm}

\begin{proof} 
First, we note that the statement directly follows when $n=4$. In this case, the 
main stratum is $\P^1 \smin \{s_1,s_2,s_3\}$ and the classes of codimension 
one strata are the boundary divisors $s_i$ which are pairwise rationally equivalent. 
The relation between the pairs of classes of codimension one strata is given by (\ref{eqn_rel}).

We prove the statement for $n>4$ by induction on $n$.

Let $\pi: \cmod{\cB} \to \cmod{\cA}$ be the map forgetting the labeled point $s_{n+1}$. Here, we use the  notations introduced in Section \ref{sec_forget}.

Let $B_d$ denote the union of $d$-dimensional (closed) strata of $\cmod{\cA}$. Let $\cmod{\cA}$ be filtered by 
\begin{eqnarray*}
\emp =B_{-1} \subset B_0 \subset \cdots  \subset B_{2n-6} = \cmod{\cA}.
\end{eqnarray*}
The forgetful map $\pi$ induces a filtration of $\cmod{\cB}$:
\begin{eqnarray*}
\emp =E_{-1} \subset E_0 \subset \cdots  \subset E_{2n-6} = \cmod{\cB}
\end{eqnarray*}
where $E_d = \pi^{-1}(B_d)$. Then, the spectral sequence obtain from this filtration 
gives us
\begin{eqnarray}
\label{eqn_d-complex}
\E^1_{p,q}      = H_{p+q} (E_p,E_{p-1})  \Longrightarrow
H_{p+q}(\cmod{\cB};\Z).
\end{eqnarray}
We prove the theorem by writing down this spectral sequence explicitly. As a first step, we calculate the homology groups  $H_{p+q} (E_p,E_{p-1})$.

From now on, we assume that the statement of the theorem holds for $\cmod{\cA}$.

\paragraph*{Step 1.} Since all strata are even dimensional, we only need to consider 
$(E_{2p},E_{2p-1})$ pairs. We can write homology of $(E_{2p},E_{2p-1})$ as a direct sum of the homology 
of its pieces:
\begin{eqnarray*}
H_{2p+q} (E_{2p},E_{2p-1}) = \bigoplus_{\BGg : \dim_\R D_{\BGg} =2p}
H_{2p+q} (\pi^{-1}(\konj{D}_{\BGg}),\pi^{-1}(Q_{\BGg})).
\end{eqnarray*}

Assume that the maximum of the dimension of the fibers 
$F_{\BGg^*}$ of $\pi: D_{\BGg^*} \to D_{\BGg}$ is  $2k$.  
Consider the following filtration of $\pi^{-1}(\konj{D}_{\BGg})$:
\begin{eqnarray*}
\emp \subset Y_0 \subset Y_1 \subset \cdots \subset Y_{2k}= \pi^{-1}(\konj{D}_{\BGg})
\end{eqnarray*}
where  $Y_i$ is the union of the strata $\konj{D}_{\BGg^*}$ of $\cmod{\cB}$ of 
dimension $2p+i$ such that $\pi(D_{\BGg^*})=D_{\BGg}$. Clearly,  
$Y_{2i}=Y_{2i+1}$ due to the absence of odd dimensional strata.

The spectral sequence of this filtration gives
\begin{eqnarray*}
\Y^1_{i,j} = H_{i+j} (Y_i,Y_{i-1} \bigcup (Y_i \cap \pi^{-1}(Q_{\BGg})))
\Longrightarrow H_{i+j} (E_{2p},E_{2p-1}).
\end{eqnarray*}
Here, $Y_i \cap \pi^{-1}(Q_{\BGg})$ contains the substrata of 
$\konj{D}_{\BGg^*}$ that maps to $B_{p-1}$ 
(i.e., substrata of codimension one or higher in  $\konj{D}_{\BGg^*}$).
Hence, we have
\begin{eqnarray*}
\Y^1_{2i,j}      = \bigoplus_{\BGg^*}  H_{2i+j} (\konj{D}_{\BGg^*},Q_{\BGg^*}), \ \ \ \ 
\Y^1_{2i+1,j}  = 0.
\end{eqnarray*}

By using Lemma \ref{lem_closed_hom} (and isomorphism between 
relative homology and homology with closed supports), we can have 
the groups $\Y_{i,j}$ as products of homology groups. We consider 
$\Y^1_{*,*}$ as a direct sum of stable and unstable pieces:

If the restriction of forgetful morphism $\pi: D_{\BGg^*} \to D_{\BGg}$ 
doesn't require stabilization, then the fibers 
$F_{\BGg^*}$ are punctured $\P^1$'s. Then, we have
\begin{eqnarray*}
\K^1_{2,j} &=& \bigoplus_{\BGg^*}  
H_j (\konj{D}_{\BGg},Q_{\BGg}) \otimes H_2 (\overline{F}_{\BGg^*},\Bs), \\
\K^1_{2,j-1} &=& \bigoplus_{\BGg^*}  
H_j (\konj{D}_{\BGg},Q_{\BGg}) \otimes H_1 (\overline{F}_{\BGg^*},\Bs).
\end{eqnarray*}
If the restriction of forgetful morphism $\pi: D_{\BGt^*} \to D_{\BGg}$ 
requires stabilization,  we have
\begin{eqnarray*}
\L^1_{i,j} = \bigoplus_{\BGt^*} 
H_j (\konj{D}_{\BGg},Q_{\BGg}) \otimes H_i (\konj{F}_{\BGt^*} ).
\end{eqnarray*}
The union of unstable fibers is $\bigsqcup_{\BGt^*} \cmod{\cB_{{v_s}}}$ 
where  $|\BV_{\Gt^*}|= |\BV_{\Gg}|+1$ and all tails supported by $v_s$ are 
distinct (i.e., there is no subset of tails that are identified). This follows from the 
fact that all other strata (corresponding to $\cB$-trees with more vertices or 
identified tails) are contained in these strata (see Proposition \ref{thm_strata}). 
Moreover, these strata are pairwise disjoint since they are uniquely determined 
by the set $\BF_{\Gt^*}(v_s)$. Hence,
\begin{eqnarray*}
\L^1_{i,j} = \bigoplus_{\BGt^* \mid |\BV_{\Gt^*}|= |\BV_{\Gg}|+1} 
H_j (\konj{D}_{\BGg},Q_{\BGg}) \otimes H_i (\cmod{\cB_{v_s}}),
\end{eqnarray*}  
and  
\begin{eqnarray}
\Y_{i,j} = \K_{i,j} \oplus \L_{i,j}.
\end{eqnarray}

Then, the differential $d_1: \Y^1_{i,j} \to \Y^1_{i-1, j}$  is zero since 
$H_{odd}(\cmod{\cB_{v_s}})$ is zero.

Finally, the differential $d_2: \Y^1_{2,j} \to \Y_{0, j+1}$ is given by the differentials
\begin{eqnarray}
\label{eqn_d2}
&\dd_*:& H_1(\overline{F}_{\BGg^*},\Bs) \to  H_0 (\overline{F}_{\BGt^*})
\end{eqnarray}
where $|\BV_{\Gt^*}|= |\BV_{\Gg}|+1$ and $|O_{\fr(v_s)}|=3$. For each pair 
of points  lying in the same component of $\konj{F}_{\BGg^*}$,  there is a generator 
in $H_1  (\overline{F}_{\BGg^*},\Bs)$ whose image under $\dd_*$ gives the difference 
of these points (see Section \ref{sec_fib_hom}). Therefore, the strata $\konj{D}_{\BGt^*_1}$, 
$\konj{D}_{\BGt^*_2}$, which are zero dimensional fibrations over $\konj{D}_{\BGg}$, are 
homologous relative to $\pi^{-1} (Q_{\BGg}$).

It is important to note that the kernel of the differential  $d_2: \Y^1_{2,j} \to \Y_{0, j+1}$ 
is trivial. This follows from the fact that the same is  true for $\dd_*$ given in (\ref{eqn_d2}) 
due to  the relations of the homology of the fibers given in (\ref{eqn_f_rel}).

By assuming the statement for the moduli space of weighted pointed curves with fewer 
labeled points, in particular for $\cmod{\cA}$, we set the factors 
$H_j (\konj{D}_{\BGg},Q_{\BGg})=\Z \cls{\BGg}$ when $j=\dim D_{\BGg}$ and zero 
otherwise  in  $\Y_{*,*}$. Hence, we obtain that $\Y_{*,*}$ are generated by the fundamental 
classes of strata. By this assumption,  the total homology of the spectral sequence $\Y$ is 
obtained by taking quotient by the image of $d_2$:  
\begin{eqnarray*}
\Y_{2i} &=&  
\left( \bigoplus_{\BGg^* \mid \dim D_{\BGg^*}=2i} \Z\ \cls{\BGg^*}  \right) / I_{2i}, \\
\Y_0 &=& 
\left( \bigoplus_{\BGt^* \mid \dim D_{\BGt^*}=2p} \Z\ \cls{\BGt^*}  \right) / I_0.
\end{eqnarray*}
The subgroups $I_{2i}$  are generated by a subset of relations 
$\cR(\BGg^*;v,f_1,f_2,f_3,f_4) \equiv 0$: 
If $\sum_{f \in \BF_{\Gg^*}(v) \smin \{s_{n+1}\}} \cB(f)  \leq 2$, they arise from 
the homology relations of 
$\cmod{\cB_{{v_s}}}$. If $\sum_{f \in \BF_{\Gg^*}(v) \smin \{s_{n+1}\}} \cB(f)  > 2$ and 
$f_i \ne s_{n+1}$, then we have again $\cR(\BGg^*;v,f_1,f_2,f_3,f_4) \equiv 0$. 
These relations are obtained by pulling back the relations from the base $\cmod{\cA}$. 
Finally, in addition to these, we have 
\begin{eqnarray}
\label{eqn_rela_rel}
\cls{\BGt^*_1} - \cls{\BGt^*_2} =0
\end{eqnarray}
because these are in the image of $d_2$ as described above.

\paragraph*{Step 2.}
The calculation in Step 1 implies that the $\E_{**}^1$ are generated by the fundamental 
classes of  the  strata. Moreover, they admit the relations that are imposed in the
statement of the theorem:  For each relation (\ref{eqn_rela_rel}) in relative homology 
$H_{*}(E_p,E_{p-1})$, there is a relation in $H_{*}(\cmod{\cB})$. We   obtain the
missing relations: $\cR(\Gg^*;f_1,f_2,f_3,s_{n+1}) \equiv 0$. This completes the 
set of relations given in the statement.

\paragraph*{Step 3.} We have a complete description of generators and relations 
in $\E^1_{*,*}$. We need to calculate the higher differentials. Since all strata of 
$\cmod{\cB}$ are even dimensional (which means $\E_{k,l}=0$ unless both $k$ and $l$ 
are even), the higher differentials 
\begin{eqnarray}
d_i: \E_{k,l} \to \E_{k-i,l+i-1}, \ \  i \geq 1 
\end{eqnarray}
must be zero. This follows from the simple fact that both $k-i$ and $l+i-1$ cannot be
zero modulo 2 when  $k,l =0\ \mod 2$.
\end{proof}

\begin{cor}
\label{thm_chow}
The Chow groups of $\cmod{\cA}$ are
\begin{eqnarray*}
A_{d} (\cmod{\cA})   &=& 
\left( \bigoplus_{\BGg \mid \ \dim_\C D_{\BGg} =d}  \Z \ \cls{\BGg} \right) / \cI_d
\end{eqnarray*}
where the subgroup of relations $\cI_d$ is generated by  $\cR(\BGg;v,f_1,f_2,f_3,f_4)$ defined 
 in (\ref{eqn_rel}). 
\end{cor}

\begin{proof}
The homology groups of complex points of the moduli space is generated by the fundamental cycles of its strata. Obviously, the strata of $\cmod{\cA}$ are algebraic cycles. Moreover, the relations between strata are given by rational equivalence in Lemma \ref{lem_relation} i.e., the Chow groups $A_i(\cmod{\cA})$ are isomorphic
to the homology groups $H_{2i} (\cmod{\cA})$.
\end{proof}



\noindent
Current Address: Max-Planck-Institute for Mathematics, Bonn, Germany \\
Email: ceyhan@mpim-bonn.mpg.de
\end{document}